\documentclass[journal]{IEEEtran}
\usepackage{cite}
\usepackage{amsmath,amssymb,amsfonts}
\usepackage{graphicx}
\usepackage{textcomp}
\usepackage{float}
\usepackage{multicol}
\usepackage{amsmath}
\usepackage{amsthm}
\usepackage{amssymb}
\usepackage{amsfonts}
\usepackage{graphicx}
\usepackage{subfigure}
\usepackage{url}
\usepackage{booktabs}
\usepackage{float}
\usepackage{color}
\usepackage{bm}

\newtheorem{Remark}{Remark}
\newtheorem{Corollary}{Corollary}

\newtheorem{Problem}{Problem}
\newenvironment{Proof}{\noindent{\em Proof:\/}}{\hfill $\Box$\par}
\newtheorem{Theorem}{Theorem}
\newtheorem{Lemma}{Lemma}

\newtheorem{Assumption}{Assumption}

\newcommand{\col}{\hbox{col}}
\newcommand{\EQQ}{\begin{eqnarray*}}
\newcommand{\ENN}{\end{eqnarray*}}
\newcommand{\EQ}{\begin{eqnarray}}
\newcommand{\EN}{\end{eqnarray}}

\ifCLASSINFOpdf
\else
\fi

\begin{document}
\title{Adaptive Output Synchronization for a Class of Uncertain Nonlinear Multi-Agent Systems over Switching Networks\footnote{.}}

\author{Jie~Huang,~\IEEEmembership{Life Fellow,~IEEE}
\thanks{This article is a slightly updated version of
\cite{Huang23}, and the work described in this article
was supported in part by the Research Grants Council of the Hong Kong Special Administrative Region under grant No. 14201621 and in part by National Natural Science Foundation of China under Project 61973260.}
\thanks{The author is with the Department of Mechanical and Automation
Engineering, The Chinese University of Hong Kong, Hong Kong (e-mail:
jhuang@mae.cuhk.edu.hk.)}}

\maketitle

\begin{abstract}
In this paper, we first study the leader-following output synchronization  problem for a class of uncertain nonlinear multi-agent systems over jointly connected  switching networks. Our approach integrates the output-based adaptive distributed observer, the conventional adaptive control technique, and the output regulation theory. Compared with the existing results,  our control law  only relies on the output of the leader instead of the state of the leader and allows the followers and the leader to have different orders.
Then,  we  further consider the rejection of a class of  bounded disturbances  with unknown bounds. Our problem includes the state consensus problem as a special case if
the followers and the leader have the same order.
\end{abstract}


\IEEEpeerreviewmaketitle

\section{Introduction}\label{Introduction}
\IEEEPARstart{T}{he} cooperative control for multi-agent systems has been one of the central control problems for nearly two decades. The cooperative control problems include  consensus or synchronization, coverage of sensor networks, coordinated motion of robotic teams, formation flying of air vehicles, distributed optimization and so on. Fundamental to all these problems is the
consensus problem.  Depending on whether or not a multi-agent system has a leader, the consensus problem can be divided into
two classes: leaderless and leader-following. The leaderless consensus problem aims to design a control law for each agent to make the states (or outputs) of all agents asymptotically synchronize to some trajectory, while the leader-following consensus problem is to design a  control law for each follower subsystem to enable
the states (or outputs) of all follower subsystems to asymptotically track a desired trajectory generated by the leader system. Regardless of the leaderless consensus problem or the leader-following consensus problem, the control law to be designed must be distributed in the sense that the control law has to satisfy some network constraints described by a communication graph.

Early results on the two consensus problems were limited to linear multi-agent systems. For example, the leaderless consensus problem was studied in \cite{Olfati1,Ren1,Seo1,Tuna1}, while the leader-following consensus problem was considered in  \cite{Hong2,Hu1,Hu3,Ni1}. Both the leaderless consensus  problem and the leader-following consensus problem were investigated in \cite{Jadbabaie1,SH3}. Then, more attention has been paid to the consensus problem of nonlinear multi-agent systems. For example, 
in \cite{Fiengo21,LiuK1,Song1,Yu1}, the consensus problem was studied for several classes of nonlinear systems 
satisfying the global Lipschitz  or the global Lipschitz-like condition.
 In \cite{DHC,LH1,SH5,WangX1}, the leader-following consensus problem was studied  via the output regulation theory and the nonlinear systems considered in \cite{DHC,LH1,SH5,WangX1} contain both disturbances and norm bounded uncertainties.   Based on the adaptive control technique, the leader-following consensus problem  was studied  for  first-order  nonlinear multi-agent systems  in \cite{Das1}, for second-order nonlinear multi-agent systems in \cite{DHC}, for multiple uncertain Euler-Lagrange systems  in \cite{CH1}, and for multiple uncertain rigid spacecraft systems in \cite{WH20}. In \cite{Das1,LG,ZhangH1},  the neural network method was used to study uncertain nonlinear multi-agent systems subject to static networks and the designed control laws  made the tracking errors uniformly ultimately bounded  for initial conditions in some prescribed compact subset. In \cite{cas23}, the leader-following consensus problem for a class of nonlinear multi-agent systems  was studied using the reinforcement learning based sliding mode control.  Nevertheless, in \cite{cas23},  the networks  are static and the disturbances are bounded functions with known bounds.
Perhaps, the most general result  on the leader-following consensus problem for nonlinear systems was given in \cite{LH2} which has the following
merits. First,
the order of each  follower  subsystem can be any positive integer and the nonlinearity does not have to satisfy the  global Lipschitz-like condition.  Second, the follower subsystems contain both constant uncertain parameters and external disturbances and  the uncertain parameters can take any constant  value, which excludes the robust control approaches such as those in \cite{DHC,LH1,WangX1}. Third,  the communication network satisfies the jointly connected condition, which is a mild condition on the communication network since it allows the network to be disconnected at any time.
Finally, compared with \cite{Das1,LG,ZhangH1}, the result in  \cite{LH2} is global and the consensus can be achieved exactly.

However, little progress has been made since the publication of \cite{LH2}.
 In this paper, we will further tackle the leader-following output synchronization problem for a class of more general nonlinear multi-agent systems than the one in \cite{LH2}
by integrating the recent results of the output-based adaptive distributed observer, the conventional adaptive control technique, and the output regulation technique.
Comparing with the state of the art of the literature such as \cite{cas23,LH2}, we will offer the following new features.
First, the leader system and all of the  follower subsystems in \cite{cas23,LH2} must have the same order while we allow the leader system and various follower subsystems to have different orders. This feature will significantly enlarge the class of leader systems.
Second,  the control laws in \cite{cas23,LH2} rely on the state of the leader system while the control law here only needs to know the output of the leader system. Thus our control law applies to the general case where the state of the leader system  is unavailable.
Third, the disturbances in \cite{LH2} are generated by a linear exosystem with a nonlinear output and the exosystem is assumed to be neutrally stable. In contrast, we will consider a class of   disturbances that can be any bounded time functions with unknown bounds, which strictly contains the disturbances in \cite{LH2}. It is also noted that
 the disturbances in \cite{cas23} are bounded functions with known bounds.


The rest of this paper is organized as follows. In Section \ref{PF}, we present some preliminaries and describe our problem. In Section \ref{MR1}, we
 solve our problem in the absence of the  disturbances. In Section \ref{MR2}, we further study our problem
with the disturbances being bounded time functions with unknown bounds.
 Finally, in Section \ref{Conclusion}, we close the paper with some concluding remarks. 
An example can be found in \cite{Huang23}.

{\bf Notation.} For any column vectors $a_i$, $i=1,...,s$, denote $\mbox{col}(a_1,...,a_s)=[a_1^T,...,a_s^T]^T$.
 We use $\sigma(t)$ to denote a piecewise constant switching signal $\sigma:[0,+\infty)\rightarrow \mathcal{P}=\{1,2,\dots,n_{0}\}$, where $n_{0}$ is a positive integer, and
$\mathcal{P}$ is called a switching index set. We assume that all switching instants $t_0=0 <t_1<t_2,\dots$ satisfy $t_{i+1}-t_{i}\geq \tau_{0}>0$ for some constant $\tau_{0}$ and all $i =  0, 1, 2, \cdots$, where $\tau_{0}$ is called the dwell time. A function $f: [t_0, \infty) \rightarrow  \mathbb{R}^n$ is said to be piecewise continuous if there exists a sequence $\{t_j,~ j = 0, 1, \ldots \}$ with a dwell time $\tau >0$ such that $f (t)$ is continuous on each time interval $[t_j, t_{j+1})$, $j = 0, 1, \ldots$.  Let
$||f||_{\infty} = \sup_{t \geq 0} f (t)$, which is called  the infinity norm of $f$.  $f$ is said to be bounded if $||f||_{\infty} $ is finite.

\section{Problem Formulation}\label{PF}
Consider a class of nonlinear multi-agent systems as follows:
\begin{equation} \label{system1}
{y}^{(r_i)}_i + f_{i}^{T}(x_{i},t)\theta_{i}  = u_i + d_i, i = 1,2,\cdots,N
\end{equation}
where $r_i \geq 1$, $u_i, y_i \in {\mathbb{R}}$, $x_i  = \col (y_i, \dot{y}_i, \cdots, {y}_i^{(r_i-1 )} )$ are the input, output and the state of the plant, respectively,
  $f_{i}:\mathbb{R}^{r_i}\times [0,+\infty)\rightarrow\mathbb{R}^{m_i}$ are known functions satisfying locally Lipschitz condition with respect to $x_{i}$ uniformly in $t$ and continuous with respect to $t$,   $\theta_{i}\in\mathbb{R}^{m_i}$ are unknown constant parameter vectors, $d_{i}: [0,+\infty) \rightarrow\mathbb{R}$ are external disturbances. In addition, there is an exosystem of the following form:
 \begin{equation}\label{exosystem1}
\dot{v}_0  = Sv_0,~~  y_0 = F v_0
\end{equation}
where $v_0 \in \mathbb{R}^n $ and $y_0 \in \mathbb{R}$ with  $S$ and $F$ two constant matrices. The output $y_0$ is to be tracked by the outputs $y_i$ of all subsystems of (\ref{system1}) in the presence of the persistent disturbances $d_i$, which can be any piecewise continuous bounded time functions with the bounds unknown, i.e., $|d_i (t)| \leq D_i$ for all $ t\geq 0$ and some unknown positive numbers $D_i$. When $N=1$, the adaptive stabilization problem of \eqref{system1} was considered in \cite{slotine}.

The system (\ref{system1}) and the exosystem (\ref{exosystem1}) together can be viewed as a multi-agent system of $(N
+1)$ agents with (\ref{exosystem1}) as the leader and the $N$
subsystems of (\ref{system1}) as $N$ followers. With respect to the plant (\ref{system1}), the exosystem (\ref{exosystem1}), and a given switching signal
$\sigma(t)$, we can define a time-varying graph
$\bar{\mathcal{G}}_{\sigma(t)}=(\bar{\mathcal{V}},\bar{\mathcal{E}}_{\sigma(t)})$\footnote{See Chapter 1 of \cite{Cai-Su-Huang-book} for
a summary of graph.} with  $\bar{\mathcal{V}}=\{0,1,\dots,N\}$ and
$\bar{\mathcal{E}}_{\sigma(t)}\subseteq \bar{\mathcal{V}}\times
\bar{\mathcal{V}}$ for all $t\geq0$, where the node $0$ is associated
with the leader system \eqref{exosystem1}  and the node $i$, $i = 1,\dots,N$,
is associated with the $i$th subsystem of system
\eqref{system1}.
For $i=1,\dots,N$, $j=0,1,\dots,N$, and $i\neq j$, $(j,i) \in
\bar{\mathcal{E}}_{\sigma(t)}$ if and only if  $u_i$ can use the
information of the $j$th agent for control at time instant $t$.
Let $\bar{\mathcal{N}}_i(t)=\{j,(j,i)\in \bar{\mathcal{E}}_{\sigma(t)}\}$
denote the neighbor set of agent $i$ at time $t$.
Let  $\mathcal{G}_{\sigma(t)}=(\mathcal{V},\mathcal{E}_\sigma(t))$  be the subgraph of $\bar{\mathcal{G}}_{\sigma(t)}$, where $\mathcal{V}=\{1,\cdots,N\}$ and $\mathcal{E}_{\sigma(t)}\subseteq\mathcal{V}\times\mathcal{V}$ is obtained from $\bar{\mathcal{E}}_{\sigma(t)}$ by removing all edges between the node $0$ and the nodes in $\mathcal{V}$.

Define an adaptive distributed dynamic compensator as follows:
\begin{subequations}\label{4.2}
\begin{align}
&{{\dot v}_i} = {S_i}{v_i} + {L_i}\sum\limits_{j\in \bar{\mathcal{N}}_i(t)}  F\left( {{v_j} - {v_i}} \right), \label{4.2a}\\
&{{\dot S}_i} = {\mu _1}\sum\limits_{j\in \bar{\mathcal{N}}_i(t)} \left( {{S_j} - {S_i}} \right), \label{4.2b}\\
&{{\dot L}_i} = {\mu _2}\sum\limits_{j \in \bar{\mathcal{N}}_i(t)} \left( {{L_j} - {L_i}} \right),\ i = 1, \ldots ,N \label{4.2c}
\end{align}
\end{subequations}
where, for $i = 1, \cdots, N$,  $v_i\in \mathbb{R}^n$, $S_i\in \mathbb{R}^{n\times n}$, $L_i\in \mathbb{R}^n$, $S_0=S$, $\mu_1, \mu_2>0$ are design parameters, and
$L_0 \in \mathbb{R}^{n}$ is a constant matrix
such that
 the following system:
\begin{equation}\label{4.6}
\dot {\bar v} = \left( {{I_N} \otimes S - {\mathcal{H}_{\sigma (t)}} \otimes \left( {L_0 F} \right)} \right)\bar v
\end{equation}
is exponentially stable, where
the matrix $\mathcal{H}_{\sigma(t)}\in \mathbb{R}^{N\times N}$ is obtained from the Laplace matrix $\mathcal{\bar L}_{\sigma(t)}$ of $\mathcal{\bar G}_{\sigma(t)}$ by removing the first row and the first column.
If there exist $L_0 \in \mathbb{R}^{n}, \mu_1>0, \mu_2>0$ such that, for any initial conditions $v_i(0),S_i(0),L_i(0),i=1,\ldots, N$ and $v_0(0)$, the solution of systems (\ref{exosystem1}) and \eqref{4.2} exists for all $t\ge 0$ and satisfies
\begin{equation}\label{4.5}
\begin{aligned}
&\mathop {\lim }\limits_{t \to \infty } ({v_i}(t) - {v_0}(t)) = 0,\\
&\mathop {\lim }\limits_{t \to \infty } ({S_i}(t) - S) = 0,\\
&\mathop {\lim }\limits_{t \to \infty } ({L_i}(t) - L_0) = 0,\ i=1,\ldots, N
\end{aligned}
\end{equation}
then,  \eqref{4.2}   is called the adaptive distributed observer for  (\ref{exosystem1}) over the graph $\mathcal{\bar G}_{\sigma(t)}$
and $L_0$ is called the observer gain matrix.

The problem of whether or not \eqref{4.2} is the adaptive distributed observer for (\ref{exosystem1}) over the graph $\mathcal{\bar G}_{\sigma(t)}$ has been extensively studied in the literature. To summarize the main result from \cite{Cai-Su-Huang-book}, let us first list the following assumptions:

\begin{Assumption}\label{ass1}
The matrix $S$ is marginally stable.
\end{Assumption}

\begin{Assumption}\label{ass2}
The pair $(F,S)$ is detectable.
\end{Assumption}

\begin{Assumption}\label{ass3}
There exists a subsequence $\{i_k\}$ of $\{i:i=0,1,2,\ldots\}$ with $t_{i_{k+1}}-t_{i_k}< \nu$ for some positive $\nu$ such that the union graph $\mathcal{\bar G}([t_{i_k},t_{i_{k+1}}))$ contains a spanning tree with the node 0 as the root.
\end{Assumption}


\begin{Assumption}\label{ass4}
$\mathcal{G}_{\sigma(t)}$ is undirected for any $t\ge 0$.
\end{Assumption}

\begin{Remark}\label{rem1}
Assumption \ref{ass1} means none of the eigenvalues of the matrix $S$ has positive real part and those eigenvalues of $S$ with zero real parts are semi-simple.
Assumption \ref{ass2} loses no generality.
 Assumption \ref{ass3} is called the jointly connected condition \cite{Ni1,SH3}. It is a mild condition on the switching graph $\mathcal{\bar G}_{\sigma(t)}$ as it allows the graph to be disconnected at every time instant. From Theorem 3.4 of \cite{Cai-Su-Huang-book}, under Assumptions \ref{ass1}-\ref{ass4}, there exists  $L_0 \in \mathbb{R}^{n}$ such that \eqref{4.6}  is exponentially stable. For the special case of Assumption  \ref{ass1} where $S$ is  neutrally stable, that is,  all the eigenvalues of the matrix $S$ have zero real parts and are semi-simple,  there exists a symmetric and positive definite matrix $R$ such that
\begin{equation}\label{4.4}
RS+S^TR=0.
\end{equation}
For this case, the matrix $L_0=\mu_0 RF^T$ with $\mu_0 >0$ is such that \eqref{4.6}  is exponentially stable.
\end{Remark}

The following result is a rephrasing of Theorem 4.7 of  \cite{Cai-Su-Huang-book}:

\begin{Lemma}\label{lem1}
  Consider systems (\ref{exosystem1})  and \eqref{4.2}. Under Assumptions \ref{ass1}-\ref{ass4}, for any $\mu_1, \mu_2>0$, and any initial conditions $v_i(0), S_i(0), L_i(0), i=1,\ldots, N$, and $v_0 (0)$, the solution of systems (\ref{exosystem1})  and \eqref{4.2} exist for all $t\ge 0$ and achieves \eqref{4.5}
exponentially.
\end{Lemma}

As a direct consequence of Lemma \ref{lem1}, we have the following result.

\begin{Corollary}\label{cor1}
  For any $k \geq 0$,
\begin{equation}\label{ui1x}
 \mathop {\lim }\limits_{t \to \infty } ( F S^{k}_i v_i - y_0^{(k)}) = 0
\end{equation}
exponentially.
\end{Corollary}

\begin{Proof}
 Since $S_i v_i - S v_0 = S_i ({v}_i - v_0)+  ({S}_i - S) v_0$, and $v_0$ and $S_i$ are bounded,  noting $y_0 = F v_0$ and  using \eqref{4.5} gives \eqref{ui1x}  for $k = 1$. For $k > 1$, first note that
\begin{equation*}
\begin{split}
& (S_i^k -S^{k}) \\
  &  = S_i^{k-1}(S_i - S) + S_i^{k-2}(S_i - S)S + \cdots  \\
  & +  S_i(S_i - S)S^{k-2} + (S_i - S)S^{k-1}.
\end{split}
\end{equation*}
Since $S_i$ are bounded, and $\mathop {\lim }\limits_{t \to \infty } (S_i -S) = 0$ exponentially,  we have
$ \mathop {\lim }\limits_{t \to \infty } (S_i^k -S^{k}) = 0$ exponentially.
Also,
\begin{equation}\label{ui1x2}
\begin{split}
 &(S_i^k v_i - S^{k} v_0) = (S_i^k v_i - S_i^k v_0 + S_i^k v_0 -S^{k} v_0)\\
 &  = S_i^k (v_i - v_0)  +  (S_i^k -S^{k}) v_0.
  \end{split}
\end{equation}
Since $v_0$ and $S_i$ are bounded, noting $y^{(k)}_0 = F S^{k} v_0$ and $ \mathop {\lim }\limits_{t \to \infty } (S_i^k -S^{k}) = 0$ exponentially, and  using \eqref{4.5} and \eqref{ui1x2} show that \eqref{ui1x}  holds for $k > 1$.
\end{Proof}

\begin{Remark}\label{rem2s}
Two special cases of Lemma \ref{lem1} are worth mentioning. First, $y_0 = v_0$, that is, the state of the leader system is available. In this case, (\ref{4.2}) reduces to the following form:
\begin{subequations}\label{4.2s}
\begin{align}
&{{\dot v}_i} = {S_i}{v_i} + \mu_v \sum\limits_{j\in \bar{\mathcal{N}}_i(t)} \left( {{v_j} - {v_i}} \right), \label{4.2as}\\
&{{\dot S}_i} = {\mu _1}\sum\limits_{j\in \bar{\mathcal{N}}_i(t)} \left( {{S_j} - {S_i}} \right), ~ i = 1, \ldots, N \label{4.2bs}
\end{align}
\end{subequations}
where $\mu_v > 0$. For this special case, Assumption \ref{ass1} can be weakened to that none of the eigenvalues of $S$ have positive real parts, Assumption \ref{ass2} is satisfied automatically, and  Assumption \ref{ass4} is not needed (see  Remark 4.4 of  \cite{Cai-Su-Huang-book}). Such a special case was considered in \cite{LH2}.
Second, the graph is static. In this case, Assumption \ref{ass3} reduces to that the graph contains a spanning tree with the node $0$ as the root, and Assumptions \ref{ass1} and \ref{ass4} are not needed (see  Theorem 4.8 of  \cite{Cai-Su-Huang-book}). Even for this special case, the problem has not been considered.
\end{Remark}

To formulate our problem, let us  first describe our control law as follows:
\begin{equation}\label{ui1}
  \begin{aligned}
 u_{i}&=h_{i}(x_{i},\zeta_{i})\\
 \dot{\zeta}_{i}&=l_{i}(x_{i},\zeta_{i},x_{j},\zeta_{j},j\in\bar{\mathcal{N}}_{i}(t)),~i=1,\cdots,N\\
  \end{aligned}
\end{equation}
where $h_{i}$ and $l_{i}$  are some nonlinear functions.

A control law of the form \eqref{ui1} is called a distributed control law since $u_{i}$ only depends on the information of its neighbors and itself. Our problem is described as follows.
\begin{Problem}\label{CORPS}
 Given the multi-agent system composed of (\ref{system1}) and (\ref{exosystem1}), bounded piecewise continuous disturbances $d_i (t)$, and  a switching graph $\bar{\mathcal{G}}_{\sigma(t)}$, design a control law of the form \eqref{ui1},  such that, for any
initial states $x_{i}(0)$, $\zeta_{i}(0)$ and $v_0(0)$,  the solution of the closed-loop system  exists  for all $t\geq0$, and satisfies
$\lim_{t\rightarrow \infty}(y^{(k)}_{i}(t)-y^{(k)}_{0}(t))=0$, $i = 1, \cdots, N$, $k = 0, \cdots, r_i-1$.
\end{Problem}

\begin{Remark}
A similar problem was considered in \cite{LH2}. However, the follower system considered in \cite{LH2} is a special case of (\ref{system1}) where $r_1 = \cdots = r_N = r$ for some $r >0$ and the leader system considered in \cite{LH2} is also a special case of (\ref{exosystem1}) where $n = r$ and $y_0 = v_0$.
In the special case where $r_1 = \cdots = r_N = n = r$, our control law actually achieves the leader-following state consensus.
It is also noted that the disturbances in \cite{LH2} are a function of trigonometric polynomials and are thus bounded. In contrast, the class of disturbances here strictly contains the disturbances in \cite{LH2} as a subclass.
It will be seen later that the class of disturbances will be treated  differently from the approach in \cite{LH2}.

\end{Remark}

\section{Adaptive Output Synchronization}\label{MR1}

In this section, we  consider the case without  disturbances.
To present our distributed control law,  let
\begin{equation}\label{pri1}
\begin{split}
p_{ri}=&~ F S^{r_i -1}_i  v_{i}\!-\! \beta_{1i}(y^{(r_i-2)}_i\!-\! F S^{r_i -2}_i  v_{i})\!-\!\cdots\\
&-\beta_{(r_i\!-\!2)i }(\dot{y}_{i}\!-\!F S_i {v}_{i}) -\beta_{(r_i\!-\!1)i }(y_{i}\!-\!F {v}_{i})
\end{split}
\end{equation}
where $\beta_{1i},\cdots,\beta_{(r_i-1)i}$ are some positive constants such that the polynomials $\lambda^{r_i-1}+\beta_{1i}\lambda^{r_i-2}+\cdots+\beta_{(r_i-2)i}\lambda+\beta_{ (r_i-1)i}=0$ are stable.
Let
\begin{equation}\label{si2}
\begin{split}
s_{i}&=y^{(r_i -1)}_i-p_{ri}.
\end{split}
\end{equation}

Now we propose our control law as follows:
\begin{equation}\label{ui3}
\begin{split}
u_{i}&=f_{i}^{T}(x_{i},t)\hat{\theta}_{i} -k_{i} s_{i}+ \dot{p}_{ri} \\
 \dot{\tilde{\theta}}_{i}&= - \Lambda_{i}^{-1}f_{i}(x_{i},t)s_{i} \\
{{\dot v}_i} &= {S_i}{v_i} + {L_i}\sum\limits_{j\in \bar{\mathcal{N}}_i(t)}  F \left( {{v_j} - {v_i}} \right)\\
{{\dot S}_i} &= {\mu _1}\sum\limits_{j\in \bar{\mathcal{N}}_i(t)} \left( {{S_j} - {S_i}} \right)\\
{{\dot L}_i} &= {\mu _2}\sum\limits_{j \in \bar{\mathcal{N}}_i(t)} \left( {{L_j} - {L_i}} \right),\ i = 1, \ldots ,N
\end{split}
\end{equation}
where $k_{i}$ are positive constants, $\hat{\theta}_{i}$ are the estimates of $\theta_i$, $\tilde{\theta}_{i} = \hat{\theta}_{i} - \theta_{i}$,
and $\Lambda_{i}\in\mathbb{R}^{m_i\times m_i}$ are symmetric and positive definite matrices.

Substituting  the control law \eqref{ui3} to the plant   gives the  closed-loop system with $d_i = 0$ as follows:
 \begin{equation}\label{system3}
\begin{split}
 \dot{x}_{li} &= {x}_{(l+1)i},~ l = 1, \cdots, r_i -1\\
    \dot{x}_{r_ii} & = f_{i}^{T}(x_{i},t)\tilde{\theta}_{i} - k_{i}s_{i}+ \dot{p}_{ri}\\
   \dot{\tilde{\theta}}_{i}&= - \Lambda_{i}^{-1}f_{i}(x_{i},t)s_{i} \\
   {{\dot v}_i} &= {S_i}{v_i} + {L_i}\sum\limits_{j\in \bar{\mathcal{N}}_i(t)}  F \left( {{v_j} - {v_i}} \right)\\
{{\dot S}_i} &= {\mu _1}\sum\limits_{j\in \bar{\mathcal{N}}_i(t)} \left( {{S_j} - {S_i}} \right)\\
{{\dot L}_i} &= {\mu _2}\sum\limits_{j \in \bar{\mathcal{N}}_i(t)} \left( {{L_j} - {L_i}} \right),\ i = 1, \ldots ,N
\end{split}
\end{equation}

Before presenting our main result, let us first establish the following result:
\begin{Lemma}
Let $e_i = y_i - y_0$. Then  \eqref{pri1} and \eqref{si2} imply
\begin{equation}\label{pri2}
e^{(r_i-1)}_i + \beta_{1i} e^{(r_i-2)}_i
\!+\!\cdots+\beta_{i (r_i\!-\!2)}\dot{e}_{i}\! + \beta_{(r_i\!-\!1)i}{e}_{i} = \bar{u}_i
\end{equation}
where
\begin{equation}\label{pri2x}
\begin{split}
&\bar{u}_i = - (y^{(r_i-1)}_0 \!-\! F S^{r_i -1}_i  v_{i}) - \beta_{1i} (y^{(r_i-2)}_0 \!-\!  F  S^{r_i -2}_i  v_{i}) \\
&\!-\!\cdots-\beta_{(r_i\!-\!2)i}  (\dot{y}_0 \!-\!F S_i {v}_{i}) -\beta_{(r_i\!-\!1)i}({y}_{0}\!-\! F {v}_{i}) + s_i
\end{split}
\end{equation}
and
\begin{equation}\label{L315}
\begin{split}
\dot{p}_{ri}
& = FS_i^{r_i}v_i \!-\! \beta_{1i}( y^{(r_i-1)}_i - FS_i^{r_i-1} v_i) \!-\!\cdots- \\
& \beta_{(r_i\!-\!2)i} ( y^{(2)}_i - FS^{2}_i v_i)
 -\beta_{(r_i\!-\!1)i} ( \dot{y}_i - FS_i v_i) \\
&+\sum_{k=0}^{r_i-1} \beta_{ki}F e_{(r_i-1-k)i}
\end{split}
\end{equation}
where
 $\beta_{0i} = 1$, $e_{0i} = e_{vi}$, and, for $k \geq 1$,
\begin{equation}\label{L315x}
\begin{split}
e_{ki} &= k S_i^{k-1}  e_{S_i}  v_i + S_i^k e_{vi}\\
e_{vi} &=  L_i \sum\limits_{j \in \bar{\mathcal{N}}_i(t)} F ({v}_{j}-{v}_{i}) \\
e_{S_i} &= \mu_1 \sum\limits_{j \in \bar{\mathcal{N}}_i(t)}   ({S}_{j}-{S}_{i}).
\end{split}
\end{equation}
\end{Lemma}
\begin{Proof}
For any $k \geq 0$,
\begin{equation}\label{L311x}
\begin{split}
& y_i^{(k)} - F S_i^{k}v_i \\
&= y_i^{(k)} - y_0^{(k)} +  y_0^{(k)} - F S_i^{k}v_i \\
&= e_i^{(k)} +   y_0^{(k)} - F S_i^{k}v_i.
\end{split}
\end{equation}
Thus,
from \eqref{pri1}, we have
\begin{equation}\label{L314}
\begin{split}
&{p}_{ri}= F S^{r_i -1}_i  v_{i}\!-\! \beta_{1i}(y^{(r_i-2)}_i\!-\! F S^{r_i -2}_i  v_{i})\!-\!\cdots\\
&-\beta_{(r_i\!-\!2)i}(\dot{y}_{i}\!-\!F S_i {v}_{i}) -\beta_{(r_i\!-\!1)i}(y_{i}\!-\!F {v}_{i})  \\
=& F S^{r_i -1}_i  v_{i} \!-\! \beta_{1i} e^{(r_i-2)}_i
\!-\!\cdots-\beta_{(r_i\!-\!2)i}\dot{e}_{i}\! -\beta_{(r_i\!-\!1)i}{e}_{i}\\
& \!-\! \beta_{1i}  (y^{(r_i-2)}_0 \!-\!  F S^{r_i -2}_i  v_{i}) \\
&\!-\!\cdots-\beta_{(r_i\!-\!2)i}  (\dot{y}_0 \!-\!F S_i {v}_{i}) -\beta_{(r_i\!-\!1)i} ({y}_{0}\!-\! F{v}_{i})
\end{split}
\end{equation}
Using \eqref{si2} and noting $y^{(r_i-1)}_i\!-\! F S^{r_i-1}_i  v_{i} = e^{(r_i-1)}_i\! + (y^{(r_i-1)}_0-\! F S^{r_i-1}_i  v_{i})$ gives \eqref{pri2}.

Next, for any $k \geq 1$,
\begin{equation}\label{L311}
\begin{split}
&\frac{d (S_i^k v_i)}{dt}
= k S_i^{k-1}  \dot{S}_i v_i + S_i^k \dot{v}_i \\
&= k S_i^{k-1} e_{S_i}  v_i + S_i^k (S_i v_i + e_{vi}) \\
&= S_i^{k+1}v_i + k S_i^{k-1}  e_{S_i}  v_i + S_i^k e_{vi} \\
&= S_i^{k+1}v_i + e_{ki}.
\end{split}
\end{equation}

Differentiating \eqref{pri1}  and using \eqref{L311}  gives
\begin{equation}
\begin{split}
\dot{p}_{ri} &= F (S^{r_i }_i  v_{i} + e_{(r_i-1)i}) \\
& \!-\! \beta_{1i}( y^{(r_i-1)}_i - FS_i^{r_i-1} v_i - F e_{(r_i - 2)i})   \\
 & \!-\!\cdots - \beta_{(r_i\!-\!2)i} ( y^{(2)}_i - FS^{2}_i v_i - F e_{1i})\\
 & -\beta_{(r_i\!-\!1)i} ( \dot{y}_i - FS_i v_i - F e_{0i})
 \end{split}
\end{equation}
which is the same as (\ref{L315}).

\end{Proof}

\begin{Remark}\label{rem3.2x}
It is interesting to note that  $\dot{p}_{ri}$ are independent of $y_{0}^{(k)}$ for $k = 0, 1, \ldots$, and thus can be used in the control law \eqref{ui3}.
\end{Remark}

Let $\xi_{1i}=e_i$, $\xi_{2i}=\dot{e}_{i}$, $\cdots$, $\xi_{(r_i-1)i}={e}_i^{(r_i-2)}$, and $\xi_{i}=\mbox{col}(\xi_{1i},\cdots,\xi_{(r_i-1)i})$. Then, \eqref{pri2}  can be put into the following state space form:
 \begin{equation}\label{system5}
\begin{split}
\dot{\xi}_{i}&=A_i \xi_{i}+ B_i \bar{u}_{i}
\end{split}
\end{equation}
where \begin{equation*}
\begin{split}
A_i=\left[\!\!\!
                    \begin{array}{cccc}
                      0 & 1 & \cdots & 0 \\
                      \vdots & \vdots & \ddots & \vdots \\
                      0 & 0 & \cdots & 1 \\
                      -\beta_{(r_i-1)i}  & -\beta_{(r_i-2)i} & \cdots & -\beta_{1i} \\
                    \end{array}
                  \!\!\!\right]\!,~B_{i}=\left[\!\!
                                         \begin{array}{c}
                                           0 \\
                                           \vdots \\
                                           0  \\
                                          1
                                         \end{array}
                                       \!\!\right].
\end{split}
\end{equation*}

Before stating our result, we need one more assumption as follows:

\begin{Assumption}\label{ass6}
$||f_i (x_i, t)|| \leq \phi_i (x_i) $ for some globally defined functions $\phi_i (x_i)$.
\end{Assumption}

Assumption \ref{ass6} is  mild since any time invariant function of $x_i$ satisfies
Assumption \ref{ass6}.

We now state the following result:
\begin{Theorem}\label{Theorem1}
Under Assumptions \ref{ass1}-\ref{ass6}, the leader-following output synchronization problem for the multi-agent system composed of \eqref{system1} and \eqref{exosystem1}
is solvable by the  control law \eqref{ui3}.
\end{Theorem}
\begin{Proof}

Let 
 \begin{equation}\label{V3}
\begin{split}
V =\frac{1}{2}\sum_{i=1}^{N}(s_{i}^{2}+\tilde{\theta}_{i}^{T}\Lambda_{i}\tilde{\theta}_{i}).
\end{split}
\end{equation}

Then,  along the solution of \eqref{system3}, we have, noting \eqref{si20},
\begin{align}
\dot{{V}} &=\sum_{i=1}^{N}(s_{i}\dot{s}_{i}+\tilde{\theta}_{i}^{T}\Lambda_{i}\dot{\tilde{\theta}}_{i} ) \notag\\
& = \sum_{i=1}^{N}\bigg(s_{i}(f_{i}^{T}(x_{i},t)\tilde{\theta}_{i}   -k_{i}s_{i} )
+\tilde{\theta}_{i}^{T}\Lambda_{i}\dot{\tilde{\theta}}_{i}
\bigg)  \notag\\
& =  - \sum_{i=1}^{N} k_{i}s^2_{i}.
\label{eq0}
\end{align}

Thus, for $i=1,\cdots,N$, $s_{i}$ and $\tilde{\theta}_{i}$ are bounded over  $[0, \infty)$.
Since ${{V}}(t)$ is lower bounded by $0$, $\lim_{t \to \infty} V(t)$ has a finite limit. We will further resort to the Barbalat's Lemma to conclude $\lim_{t \to \infty} \dot{V} (t) = 0$.  For this purpose, we need to show
\begin{equation}\label{ddotV2}
\begin{split}
\ddot{V}=& - \sum_{i=1}^{N} 2 k_i s_{i}\dot{s}_{i}
\end{split}
\end{equation}
is bounded. Note that \eqref{pri2}  can be put into the form \eqref{system5}.
From \eqref{pri2x},
$\bar{u}_{i} $  is bounded
since $s_{i}$, $S_i$ and $v_{i}$ are  bounded.  Since $A_i$ are all Hurwitz, both $\xi_i$ and $\dot{\xi}_i$ are bounded, which imply
 $e_{i}^{(r_i-1)},e_{i}^{(r_i-2)}\cdots, e_{i}$ are all bounded, and hence  $x_{i}$ are bounded since $y_{0}^{(r_i-1)},y_{0}^{(r_i-2)}\cdots, y_{0}$ are. By Assumption \ref{ass6}, $f_i (x_i, t)$ are bounded since $x_{i}$ are.
 Differentiating \eqref{si2} and using the first equation of the control law \eqref{ui3}  gives
\begin{equation}\label{si20}
\begin{split}
\dot{s}_{i}&=y^{(r_i)}_i-\dot{p}_{ri} \\
&= u_i  - f_{i}^{T}(x_{i},t) {\theta}_{i} -\dot{p}_{ri} \\
& =f_{i}^{T}(x_{i},t)\tilde{\theta}_{i} - k_{i}s_{i}.
\end{split}
\end{equation}
which implies  $\dot{s}_{i}$ are also bounded.
Since $s_{i}$ and $\dot{s}_{i}$  are both bounded,  it follows from \eqref{ddotV2} that $\ddot{V}$ is bounded over $[0, \infty)$.

Thus, by the Barbalat's Lemma,  $$\lim_{t\rightarrow \infty}\dot{V}(t)=0,$$ and hence  $\lim_{t\rightarrow \infty}s_{i}(t)=0$ for $i=1,\cdots,N$. By Corollary \ref{cor1} and the fact that $\lim_{t\rightarrow \infty}s_{i}(t)=0$,  we have $\lim_{t\rightarrow \infty} \bar{u}_{i}(t)=0$ for $i=1,\cdots,N$.
Since $A_i$ are Hurwitz, from \eqref{system5},
we have $\lim_{t\rightarrow \infty}\xi_{i}(t)=0$ and $\lim_{t\rightarrow\infty} \dot{\xi}_{i}(t)=0$, which imply  $\lim_{t\rightarrow \infty} e^{(k)}_i (t) =0$ for $k=0,\ldots,r_i-1$ and $i=1,\ldots,N$.
\end{Proof}

By Remark \ref{rem2s}, for the special case where the graph is static, we have the following corollary of Theorem \ref{Theorem1}.

\begin{Corollary}
Under Assumptions \ref{ass1}, \ref{ass2},  \ref{ass6} and the assumption that  the graph is static and contains a spanning tree with the node $0$ as the root, the leader-following output synchronization  problem for the multi-agent system composed of \eqref{system1} and \eqref{exosystem1} is solvable by the  control law \eqref{ui3}.
\end{Corollary}

\section{Disturbance Rejection}\label{MR2}

In this section, we further consider the  disturbance rejection problem.
Instead of \eqref{ui3}, we modify the control law to the following form:
\begin{equation}\label{ui32}
\begin{split}
u_{i}&=f_{i}^{T}(x_{i},t)\hat{\theta}_{i}-sgn(s_i) \hat{D}_i -k_{i}s_{i}+ \dot{p}_{ri} \\
\dot{\tilde{\theta}}_{i}&= - \Lambda_{i}^{-1}f_{i}(x_{i},t)s_{i}\\
\dot{\tilde{D}}_{i}&= sgn (s_i) s_{i}\\
 {{\dot v}_i} &= {S_i}{v_i} + {L_i}\sum\limits_{j\in \bar{\mathcal{N}}_i(t)}  F \left( {{v_j} - {v_i}} \right)\\
{{\dot S}_i} &= {\mu _1}\sum\limits_{j\in \bar{\mathcal{N}}_i(t)} \left( {{S_j} - {S_i}} \right)\\
{{\dot L}_i} &= {\mu _2}\sum\limits_{j \in \bar{\mathcal{N}}_i(t)} \left( {{L_j} - {L_i}} \right),\ i = 1, \ldots ,N
\end{split}
\end{equation}
where $\hat{D}_i$ are the estimates of the upper bounds  $D_i$ of $d_i$, $\tilde{D}_{i} = \hat{D}_i - D_i$,  and, for any scalar $x$,  the function $sgn (\cdot)$ is defined as follows:
\begin{align}
sgn(x)&=\begin{cases}
1,\quad &x> 0 \\
0,\quad &x= 0 \\
-1,\quad &x<0.
\end{cases} \label{sign_fun0}
\end{align}

Under the control law \eqref{ui32}, the closed-loop system can be put as follows:
 \begin{equation}\label{system32}
\begin{split}
 \dot{s}_{i} & = f_{i}^{T}(x_{i},t)\tilde{\theta}_{i}+{d}_{i} -sgn(s_i) \hat{D}_i - k_{i}s_{i}\\
    \dot{\tilde{\theta}}_{i}&= - \Lambda_{i}^{-1}f_{i}(x_{i},t)s_{i}\\
      \dot{\tilde{D}}_{i}&= sgn (s_i) s_{i}\\
    {{\dot v}_i} &= {S_i}{v_i} + {L_i}\sum\limits_{j\in \bar{\mathcal{N}}_i(t)}  F \left( {{v_j} - {v_i}} \right)\\
{{\dot S}_i} &= {\mu _1}\sum\limits_{j\in \bar{\mathcal{N}}_i(t)} \left( {{S_j} - {S_i}} \right)\\
{{\dot L}_i} &= {\mu _2}\sum\limits_{j \in \bar{\mathcal{N}}_i(t)} \left( {{L_j} - {L_i}} \right),\ i = 1, \ldots ,N.
 \end{split}
\end{equation}

It is noted that the right hand side of the closed-loop system \eqref{system32} is discontinuous in $s_i$. Thus, the solution of the closed-loop system \eqref{system32}  must  be defined in the Filipov sense \cite{cortes,shevitz1994lyapunov}.  Put the first three equations of the closed-loop system \eqref{system32} to the following compact form:
\begin{align} \label{eqc}
\dot{x}_{c}  = f_c (x_c, t).
\end{align}
where $x_c = \col (s_1, \cdots, s_N, \tilde{\theta}_{1}, \cdots, \tilde{\theta}_{N}, \tilde{D}_{1}, \cdots, \tilde{D}_{N})$.  Then,  the Filipov solution of \eqref{eqc} satisfies, for almost all $t \geq 0$,
\begin{align} \label{eqF}
 \dot{x}_{c}  \in K [f_c] (x_c, t)
\end{align}
where $K [f_c] (x_c, t)$ is the Filipov set  of $f_c (x_c, t)$ \cite{cortes,shevitz1994lyapunov}. It is known that, for any scalar $x$,
the Filipov set  of $sgn (x)$  denoted by $K[sgn](x)$ is as follows \cite{shevitz1994lyapunov}:
\begin{align}
K[sgn] (x)&=\begin{cases}
1,\quad &x> 0 \\
[-1,1],\quad &x= 0 \\
-1,\quad &x<0.
\end{cases} \label{sign_fun}
\end{align}
Thus
\begin{align} \label{id1}
x K[sgn] (x)&= |x|.
\end{align}

We have the following result:
\begin{Theorem}\label{Theorem3}
Under Assumptions \ref{ass1}-\ref{ass6},   the leader-following output synchronization with disturbance rejection problem for the multi-agent system composed of \eqref{system1} and \eqref{exosystem1} is solvable by the distributed control law \eqref{ui32}.
\end{Theorem}

\begin{Proof}
Let 
 \begin{equation}\label{V32}
\begin{split}
V =\frac{1}{2}\sum_{i=1}^{N}(s_{i}^{2}+\tilde{\theta}_{i}^{T}\Lambda_{i}\tilde{\theta}_{i} + \tilde{D}^2_i)
\end{split}
\end{equation}
whose gradient is $$\partial V= [s_1, \cdots, s_N, \tilde{\theta}_{1}^{T}\Lambda_{1}, \cdots, \tilde{\theta}_{N}^{T}\Lambda_{N},   \tilde{D}_1, \cdots, \tilde{D}_N].$$
By Theorem 2.2 of \cite{shevitz1994lyapunov}, $\dot{V}$ exists almost everywhere (a.e.), and $\dot{V}\in ^{a.e.}\dot{\tilde{V}}$, where
\begin{align}
	\dot{\tilde{V}}&=  \begin{array}{cc} \bigcap & \psi K [f] ({x}_c, t) \\
                                      \footnotesize{ \psi \in    \partial V}  & ~\end{array}     \notag\\
&=\sum_{i=1}^{N}(s_{i}\dot{s}_{i}+\tilde{\theta}_{i}^{T}\Lambda_{i}\dot{\tilde{\theta}}_{i} + \tilde{D}_i \dot{\tilde{D}}_i) \notag\\
=& \sum_{i=1}^{N}\bigg(s_{i}(f_{i}^{T}(x_{i},t)\tilde{\theta}_{i}+{d}_{i}-K [sgn] (s_i) \hat{D}_i-k_{i}s_{i})\notag \\
&-\tilde{\theta}_{i}^{T}f_{i}(x_{i},t)s_{i} + \tilde{D}_i K[sgn](s_i) s_i \bigg) \notag \\
= & \sum_{i=1}^{N}\bigg(s_{i}{d}_{i}-s_i K[sgn](s_i){D}_i-k_{i}s^2_{i} \bigg).
 \label{eq02}
\end{align}
Using (\ref{id1}) in \eqref{eq02} gives
\begin{align}
\dot{\tilde{V}}& = \sum_{i=1}^{N} (s_i d_i - |s_i| D_i -  k_{i} s_{i}^{2}).
\end{align}
Thus, $\dot{\tilde{V}} = \dot{{V}}$. Noting $(s_i d_i -|s_i| D_i) \leq 0$ gives
\begin{align}
\dot{{V}}& \leq  - \sum_{i=1}^{N} k_{i} s_{i}^{2}.
\end{align}

Thus,  $ \lim_{t \to \infty} V (t)$ exists and is finite. Hence, for $i=1,\cdots,N$, $s_{i}$, $\tilde{\theta}_{i}$,  and $\tilde{D}_{i}$  are bounded over $[0, \infty)$.
Let \begin{align*}
	W(t)=\int_{0}^{t} \sum_{i=1}^{N} k_{i} s_{i}^{2}   (\tau)d\tau.
\end{align*}
Then, for all $t \geq 0$,
\begin{align*}
	W(t) \leq -  \int_{0}^{t} \dot{{V}}(\tau)d\tau=-{{V}}(t)+{{V}}(0).
\end{align*}
Since ${{V}}(t)$ is lower bounded, $\lim_{t \to \infty}W(t)$ has a finite limit.  Similar to the proof of Theorem \ref{Theorem1},  we will resort to the  Barbalat's Lemma again. For this purpose, we need to show
\begin{equation}\label{ddotW2}
\begin{split}
\ddot{W}=&  \sum_{i=1}^{N} 2 k_i s_{i}\dot{s}_{i}
\end{split}
\end{equation}
is bounded. Consider  \eqref{system5} again.
From \eqref{pri2x},
$\bar{u}_{i} $  is bounded
since $s_{i}$, $S_i$ and $v_{i}$ are  bounded.  Since $A_i$ are all Hurwitz, both $\xi_i$ and $\dot{\xi}_i$ are bounded, which imply
 $e_{i}^{(r_i-1)},e_{i}^{(r_i-2)}\cdots, e_{i}$ are all bounded, and hence  $x_{i}$ are bounded since $y_{0}^{(r_i-1)},y_{0}^{(r_i-2)}\cdots, y_{0}$ are. By Assumption \ref{ass6}, $f_i (x_i, t)$ are bounded since $x_{i}$ are. Also, $d_i$ are bounded by assumption.
From the first equation of \eqref{system32}, we have
\begin{equation*}\label{si20x}
s_{i} \dot{s}_{i}=s_{i}f_{i}^{T}(x_{i},t)\tilde{\theta}_{i} +s_{i}{d}_{i} -s_{i}sgn(s_i) \hat{D}_i - k_{i}s^2_{i},
\end{equation*}
which implies $\ddot{W}$ is  also bounded.
However, since $d_i$ and hence $\ddot{W} $ may  be discontinuous at infinitely many time instances, we cannot  invoke
the Barbalat's Lemma.  Nevertheless, we have shown that $W (t)$ satisfies the following three conditions:
\begin{itemize}
\item [1)]
$\lim_{t \to \infty} W (t)$ exists;
\item [2)]
$W (t)$ is twice differentiable on each time interval $[t_j,t_{j+1})$ satisfying $t_{j+1}-t_j\geq \tau >0$;
\item [3)]
$\ddot{W}_i(t)$ is bounded over $[0, \infty)$ in the sense that there exists a finite positive constant $K$ such that
\begin{equation}
\sup_{t_j\leq t \leq t_{j+1}, j=0,1,2,\cdots} |\ddot{W}_i(t)|\leq K.
\end{equation}
\end{itemize}
Thus, by the generalized Barbalat's Lemma as can be found in Corollary 2.5 of \cite{Cai-Su-Huang-book} to conclude $\lim_{t \to \infty} \dot{W} (t)=0$, and hence $\lim_{t\rightarrow \infty}s_{i}(t)=0$ for $i=1,\cdots,N$. By Corollary \ref{cor1} and the fact that $\lim_{t\rightarrow \infty}s_{i}(t)=0$,  we have $\lim_{t\rightarrow \infty} \bar{u}_{i}(t)=0$ for $i=1,\cdots,N$.
From \eqref{system5},
we have $\lim_{t\rightarrow \infty}\xi_{i}(t)=0$ and $\lim_{t\rightarrow \infty} \dot{\xi}_{i}(t)=0$, which implies   $\lim_{t\rightarrow \infty} e^{(k)}_i (t)=0$ for $k=0,\ldots,r_i-1$ and $i=1,\ldots,N$.
\end{Proof}

By Remark \ref{rem2s}, for the special case where the graph is static, we have the following corollary of Theorem \ref{Theorem3}.

\begin{Corollary}
Under Assumptions \ref{ass1}, \ref{ass2}, \ref{ass6}, and the assumption that  the graph is static and contains a spanning tree with the node $0$ as the root, the leader-following output synchronization with disturbance rejection problem for the multi-agent system composed of \eqref{system1} and \eqref{exosystem1}  is solvable by the distributed control law \eqref{ui32}.
\end{Corollary}

\begin{Remark}
If $d_i$ are discontinuous at only finitely many time instances, then it suffices to invoke Barbalat's lemma to conclude the proof.
\end{Remark}

\begin{Remark}
Compared with  \cite{cas23,LH2}, this paper offers the following new features.
First, the leader system and all of the  follower subsystems in \cite{cas23,LH2} must have the same order while we allow the leader system and various follower subsystems to have different orders. This feature will significantly enlarge the class of leader systems.
Second,  the control laws in \cite{cas23,LH2} rely on the state of the leader system while the control law here only needs to know the output of the leader system. Thus our control law applies to the general case where the state of the leader system  is unavailable.
Third,  the disturbances in \cite{LH2} are generated by a neutrally stable linear exosystem with a nonlinear output while our disturbance can be any bounded function.
\end{Remark}

\section{Conclusion}\label{Conclusion}
In this paper, we have studied the leader-following output synchronization problem for a class of higher-order nonlinear multi-agent systems subject to both constant parameter uncertainties and external disturbances over jointly connected switching networks.
The class of disturbances includes any bounded piecewise continuous time function with unknown bounds.
We have solved our problem by integrating the output-based adaptive distributed observer, the conventional adaptive control technique, and the adaptive disturbance rejection technique.
Our problem includes the state consensus problem as a special case if
the followers and the leader have the same order. The disturbances considered in this paper satisfy the matching condition.

\end{document}